\documentclass{siamart1116}

\usepackage{algorithm}
\usepackage{amscd}
\usepackage{amsmath}
\usepackage{amsopn}
\usepackage{amssymb}
\usepackage[noend]{algpseudocode}
\usepackage{bm}
\usepackage{cancel}
\usepackage{cleveref}
\usepackage{color}
\usepackage{enumerate}
\usepackage{float}
\usepackage{graphicx}
\usepackage[utf8]{inputenc}
\usepackage{mathrsfs}
\usepackage{mathtools}
\usepackage{multirow}
\usepackage{nccmath}
\usepackage{subcaption}
\usepackage{wrapfig}

\newcommand{\abs}[1]{\left\lvert#1\right\rvert}

\newcommand{\N}{\numberset{N}}

\newcommand{\numberset}{\mathbb}

\renewcommand{\ge}{\geqslant}
\renewcommand{\le}{\leqslant}

\DeclareSymbolFont{eulargesymbols}{U}{zeuex}{m}{n}
\DeclareMathSymbol{\intop}{\mathop}{eulargesymbols}{"52}
\DeclareMathSymbol{\ointop}{\mathop}{eulargesymbols}{"49}

\makeatletter
\newcommand{\srcsize}{\@setfontsize{\srcsize}{5pt}{5pt}}
\newcommand{\eqnum}{\refstepcounter{equation}\textup{\tagform@{\theequation}}}
\makeatother

\DeclareMathOperator{\height}{height}

\newsiamthm{Case}{Case}
\newsiamthm{Corollary}{Corollary}
\newsiamthm{Definition}{Definition}
\newsiamthm{Lemma}{Lemma}
\newsiamthm{Proposition}{Proposition}
\newsiamthm{Remark}{Remark}
\newsiamthm{Theorem}{Theorem}
\newsiamthm{Observation}{Observation}

\numberwithin{theorem}{section}

\newcommand{\TheTitle}{The $\beta$ maps: strong clustering and distribution results on the complex unit circle}
\newcommand{\TheAuthors}{A.J.A. Schiavoni Piazza, D. Meadon, S. Serra-Capizzano}

\headers{da inserire in caso}{\TheAuthors}

\title{\TheTitle}

\author{Alec J.A. Schiavoni Piazza \thanks{Department of Science and High Technology, University of Insubria, via Valleggio, 11 22100 Como, Italy (\email{ajapiazza@studenti.uninsubria.it}.} 
\and David Meadon \thanks{Division of Scientific Computing, Department of Information Technology, University of Uppsala, - Lägerhyddsv. 1, Box 337, SE-751 05, Uppsala, Sweden (\email{david.meadon@uu.se}).}
\and Stefano Serra-Capizzano \thanks{Department of Science and High Technology, University of Insubria, via Valleggio, 11 22100 Como, Italy; Division of Scientific Computing, Department of Information Technology, University of Uppsala, - Lägerhyddsv. 1, Box 337, SE-751 05, Uppsala, Sweden (\email{s.serracapizzano@uninsubria.it}).}}

\ifpdf\hypersetup{
pdftitle={\TheTitle},
pdfauthor={\TheAuthors}}
\fi

\usepackage{fancyhdr}
\begin{document}

\maketitle

\begin{abstract}
In the current work, we study the eigenvalue distribution results of a class of non-normal matrix-sequences which may be viewed as a low rank perturbation, depending on a parameter $\beta>1$, of the basic Toeplitz matrix-sequence $\{T_n(e^{\mathbf{i}\theta})\}_{n\in\N}$, $\mathbf{i}^2=-1$. The latter of which has obviously all eigenvalues equal to zero for any matrix order $n$, while for the matrix-sequence under consideration we will show a strong clustering on the complex unit circle. A detailed discussion on the outliers is also provided.
The problem appears mathematically innocent, but it is indeed quite challenging since all the classical machinery for deducing the eigenvalue clustering does not cover the considered case. In the derivations, we resort to a trick used for the spectral analysis of the Google matrix plus several tools from complex analysis. We only mention that the problem is not an academic curiosity and in fact stems from problems in dynamical systems and number theory. Additionally, we also provide numerical experiments in high precision, a distribution analysis in the Weyl sense concerning both eigenvalues and singular values is given, and more results are sketched for the limit case of $\beta=1$.

\end{abstract}



\section{Introduction}
Consider a class of Toeplitz matrix-sequences $\{T_n(e^{\mathbf{i}\theta})\}_{n\in\N}$, $\mathbf{i}^2=-1$, which have a low rank perturbation depending on a parameter $\beta>1$. We will show in this paper that for this specific class of non-normal matrix-sequences that whilst the basic Toeplitz sequence has clearly all eigenvalues equal to $0$, the family of perturbed $\beta$-matrix sequences shares a strong clustering on the complex unit circle, i.e the the range of the generating function $e^{\mathbf{i}\theta}$. For $\beta\geq 2$ no outliers show up, while for $\beta \in (1,2)$ there are only two outliers, which are both real, positive, and have a finite limit equal to $\beta-1$ and $(\beta-1)^{-1}$, respectively. We note that the considered problem is not simply an academic curiosity but does indeed stem from problems in dynamical systems and number theory (see e.g. \cite{Parry,Verger-G}).
As the usual techniques for deducing eigenvalue clustering results are not easy to apply in the current setting, the problem becomes quite challenging; Thus, at most, we may hope for a weak clustering: for the notion of strong and weak clustering see Definition \ref{def:cluster}, while for ad hoc results see \cite{appl-maj II,Barbarino-SSC,Golinskii-SSC,Lub,lothar,grande Alec-1,taud2,appl-maj I,TilliNonH,Tyrty-cl1,Tyrty-cl2} and the tools in the books \cite{BoGr00,BoGr05,BS,book GLT I,book GLT II}. As such, we will resort to a careful analysis of the characteristic polynomials as well as several tools from complex analysis in the ensuing derivations. Below, we report a concise account on $\beta$ maps from which it can be seen that our main results improve substantially those in the relevant literature.

Taking $\beta>1$ a fixed parameter, the corresponding $\beta$ map can be defined. Indeed, the latter is a function from $[0,1]$ to $[0,1]$ which multiplies by $\beta$ the input and takes the fractional part. In the model, the parameter $\beta$ is chosen in such a way that the orbit of $1$ is periodic and the period of the orbit is exactly the size $n$ of the matrix under consideration \cite{Parry}. 

In \cite[Theorem 3.2]{Verger-G} the author shows the weak clustering of the eigenvalues. Here we report selected links between the different notations and results.
\begin{itemize}
\item We consider only Parry numbers \cite{Parry} and hence the theorem can be applied only in the case where the matrix can be written;
\item the roots of the polynomial $n^*(\beta)$ are exactly the eigenvalues of our matrix;
\item $d_P$ is the order of the square matrix which in our setting is $n$; 
\item the term $\log(\beta)$ is always positive in our setting since $\beta>1$;
\item in \cite[Remark 3.3]{Verger-G} we see that the result implies that the proportion of the the eigenvalues outside a $\epsilon$ fattening of the complex unit circle is bounded by $O(\log(n)/n)$. 
\end{itemize} 

The main finding in the present work is to improve the estimates in \cite{Verger-G}. More in detail, we reduce $O(\log(n))$ to zero for $\beta\ge 2$ and to $2$ for $\beta \in (1,2)$, so passing from the weak clustering proven in \cite{Verger-G} to the strong clustering, according to Definition \ref{def:cluster}. Furthermore, in the case where $\beta \in (1,2)$, the two outliers show a favorable smooth behavior, since they both have limits as $n$ tends to infinity, equal to $\beta-1$ and $(\beta-1)^{-1}$, respectively. Additional eigenvalue and singular value distribution results are given at the end of the paper.

The current work is organized as follows. In Section \ref{notation-def} we introduce the basic notations and definitions. Section \ref{discussion} contains a concise discussion concerning the tools for proving eigenvalue clustering results. Section \ref{main} is devoted to the main derivations and findings, while in Section \ref{end} we discuss numerical experiments and present further distribution results, possible extensions, as well as open problems.

\section{Notations and Definitions}\label{notation-def}

First, we introduce the definitions, notation, and mathematical objects involved in the following sections.  
\begin{Definition}\label{Toeplitz-sequence}
	Take $f\in L^1(-\pi,\pi)$ and let $T_n(f)$ be the Toeplitz matrix generated by $f$ i.e.
$\left(T_n(f)\right)_{s,t}=\hat f_{s-t}$, $s,t=1,\ldots,n$, with $f$ denoted as the generating function of $\{T_n(f)\}_{n\in\N}$ and
with $\hat f_k$ being the $k$-th Fourier coefficient of $f$, that is
\[
\hat f_k={1\over 2\pi}\int_{-\pi}^{\pi} f(\theta)\ e^{-\mathbf{i}k \theta}\, \mathrm{d}\theta,\ \ \ \mathbf{i}^2=-1, \ \ k\in \mathbb Z.
\]

\end{Definition}

If $f$ is real-valued then several spectral properties (localization, extremal behaviour, collective distribution) are known
(see \cite{book GLT I,book GLT II} and references therein) and $f$ is also the spectral symbol of $\{T_n(f)\}_{n\in\N}$ in the Weyl sense \cite{BS,book GLT I,book GLT II,TilliL1,TyrtyZ}.
If $f$ is complex-valued, then the similar information is transferred to the singular values \cite{book GLT I,book GLT II}, while the eigenvalues may exhibit erratic behaviour \cite{appl-maj I} in some cases and regular behaviour in others \cite{TilliNonH}.

\begin{Definition}\label{Beta-mat-seq}
	Take $\beta\in(1,+\infty)$. Given a fixed $n\in\N$, let $v\in \mathbb{C}^n$ be the vector defined by
	$v_j = \beta^{-j}$, $j=1,\dots,n$, and $e_j$ be the j-th vector of the canonical basis of $\mathbb{C}^n$, and $e=\sum_{j=1}^{n}e_j$, where all vectors are considered as column vectors. 
	We then define the $\beta$ matrix of order $n$ as follows
		\begin{equation}\label{Bn compact}
		B_n=T_n(e^{\mathbf{i}\theta})+(v-e_1)e^{T}.
		\end{equation}
	The $\beta$ matrix of order $n$ can be written explicitly as
	 	\begin{equation}\label{Bn explicit}
	 	B_n=\begin{pmatrix}
	 	\beta^{-1}-1&\beta^{-1}-1&\beta^{-1}-1&\beta^{-1}-1&\cdots&\beta^{-1}-1\\
	 	\beta^{-2}+1&\beta^{-2}&\beta^{-2}&\beta^{-2}&\cdots&\beta^{-2}\\
	 	\beta^{-3}&\beta^{-3}+1&\beta^{-3}&\beta^{-3}&\cdots&\beta^{-3}\\
	 	\beta^{-4}&\beta^{-4}&\beta^{-4}+1&\beta^{-4}&\ddots&\beta^{-4}\\
	 	\vdots&\ddots&\ddots&\ddots&\ddots&\vdots\\
	 	\beta^{-n}&\cdots&\cdots&\beta^{-n}&\beta^{-n}+1&\beta^{-n}
	 	\end{pmatrix}.
	 	\end{equation}
	We are interested in the study of the spectral properties of the $\beta$ matrix-sequence $\{B_n\}_{n\in\N}$, specifically in a clustering result of the eigenvalues.
\end{Definition}

\begin{Definition}\label{def:cluster}
	A matrix-sequence $\{X_n\}_{n\in\N}$, where $X_n$ is a square matrix of size $n$, is strongly clustered at $s\in\mathbb{C}$ 
(in the eigenvalue sense) if, for any $\epsilon>0$, the number of eigenvalues of $X_n$, whose distance from $s$ is larger than $\epsilon$, is of order $O(1)$. An analogous definition can be given for the property of being strongly clustered at a closed subset $S$ of $\mathbb{C}$, 
with the distance from a point is replaced by the usual distance from a subset of a metric space.
If instead of $O(1)$ we had written $o(n)$, then the clustering would be called weak.
\end{Definition}

\section{Clustering Tools}\label{discussion}

In this brief section we collect and discuss tools for deducing eigenvalue clustering results.

For instance, when Hermitian matrix-sequences are considered, the assumptions that $\{Y_n\}_{n\in\N}$ is weakly clustered at $S$ and $\{W_n\}_{n\in\N}$ is weakly clustered at $0$ imply that $\{X_n\}_{n\in\N}$, with $X_n=Y_n+W_n$, is weakly clustered at $S$. However this result (see e.g. \cite[Exercise 5.3]{book GLT I}) requires that all the involved matrices are Hermitian, and we are clearly not in this setting.
The notion of quasi-Hermitian matrix-sequences allows us to enlarge the scope of the previous result, 
by requiring that some Schatten $p$ norms \cite{bhatia}, $p\in [1,\infty]$, of the difference $X_n-X_n^*$, with $X_n^*$ being the transpose conjugate of $X_n$, do not grow too fast 
(e.g. at most $o(n^{1/p})$, see \cite{Barbarino-SSC,Golinskii-SSC}). 
Unfortunately, there is no such value of $p$ such that the previous hypothesis holds for the considered matrix-sequences in this paper.

On the other hand, even if the matrix-sequence is close to normal, for instance a rank-one infinitesimal perturbation of a normal matrix-sequence, indeed the behavior of the resulting distribution and of the related clustering set can be discontinuous, see \cite[p. 84]{appl-maj I}. 

The previous short discussion informs us that the Toeplitz and GLT machinery may not be effective in the current setting, even if the tools introduced in \cite{TilliNonH} could still be used for deducing weak clustering of the eigenvalues.   

However, we are interested in the strong clustering and so, as done in the study of the Jordan canonical form of the Google matrix \cite{google}, we return to the basic idea of writing the eigenvalue-eigenvector relation
\[
	B_n x = [T_n(e^{\mathbf{i}\theta})+(v-e_1)e^{T}] x= \lambda x,
\]
of the matrix $B_n$ which leads to a two-term recurrence relation for the entries of the eigenvector $x$. By imposing $e^T x=1$ (the same trick as done with the Google matrix, $e$ being the vector of all ones), we arrive to an expression for all the entries of the eigenvector $x$ as a function of $\lambda$, in which the characteristic polynomial of $B_n$ is encoded. 
In the following section we will use this characteristic polynomial of $B_n$ to further continue the analysis of the eigenvalue clustering.

\section{Main Results}\label{main}

\begin{Theorem}\label{thm-strong-cluster}
	Let $\{B_n\}_{n\in\N}$ be the $\beta$ matrix-sequence given in Definition \ref{Beta-mat-seq} and let $S^1$ be the unit circle in $\mathbb{C}$. Then, the sequence $\{B_n\}_{n\in\N}$ is strongly clustered (in the eigenvalue sense) at $S^1$. In addition, if $\beta\in\left[2,+\infty\right)$ then for every given $\epsilon>0$, and $n$ large enough, no eigenvalues of $B_n$ lies at a distance greater than $\epsilon$ from $S^1$. In the other case, namely $\beta\in\left(1,2\right)$, the same condition holds for all but two eigenvalues of $B_n$. The two outliers, for $n$ large enough, are real, positive, and converge to $(\beta-1)^{-1}$ and $\beta-1$, respectively.  
\end{Theorem}
In order to prove Theorem \ref{thm-strong-cluster}, we will first compute the characteristic polynomial $p_n$ of $B_n$, for every $n$, and then proceed with an analysis of the sequence of polynomials $\{p_n\}_{n\in \N}$, taking the limit when possible and deriving information on the distribution of the roots. The resulting analysis, and an application of Hurwitz' Theorem (see Theorem \ref{Hurwitz}), will then be used to prove Theorem \ref{thm-strong-cluster}.

\begin{Lemma}\label{aux-determinant}
	Let $t\in\mathbb{C}$ and $n\geq1$. We take the auxiliary matrix $M_n = -I_n + t(T_n(e^{-\mathbf{i}\theta})-e_ne^{T})$, that is, the following matrix
		\begin{equation*}
		M_n=\begin{pmatrix}
		-1&t&0&0&\cdots&0\\
		0&-1&t&0&\cdots&0\\
		0&0&-1&t&\ddots&0\\
		\vdots&\cdots&\ddots&\ddots&\ddots&\vdots\\
		0&\cdots&\cdots&0&-1&t\\
		-t&\cdots&\cdots&\cdots&-t&-1-t
		\end{pmatrix}.
		\end{equation*}
	The determinant of $M_n$ is given by $(-1)^{n}\sum_{i=0}^nt^i$.
\end{Lemma}  
\begin{proof}
	We prove the lemma by induction, using the Laplace expansion along the first column, exploiting the fact that the cofactor of order $(1,1)$ of $M_n$ is $M_{n-1}$ and the cofactor of order $(n,1)$ is a lower triangular matrix with all diagonal entries equal to $t$.
	Base case:
	\begin{equation*}
		{\rm det}(M_1)={\rm det}\begin{pmatrix}
		-1-t\\
		\end{pmatrix}=(-1-t)=(-1)^1\sum_{j=0}^{1}t^j.
	\end{equation*}
	Next, taking into account the previous observations, we proceed with the inductive step to complete the proof. More precisely, we have
	\begin{eqnarray*}
	{\rm det}(M_n)=\sum_{i=1}^{n}(-1)^{i+1}m_{i,1}{\rm det}((M_n)_{i,1})=\\
	=(-1)^2(-1){\rm det}(M_{n-1})+(-1)^{n+1}(-t)t^{n-1}=\\
	=(-1)(-1)^{n-1}\sum_{i=0}^{n-1}t^i+(-1)^nt^n=(-1)^n\sum_{i=0}^{n}t^i,
	\end{eqnarray*}
	where $m_{i,j}$ is the element in position $(i,j)$ of $M_n$ and $(M_n)_{i,j}$ is the cofactor of order $(i,j)$ of $M_n$. This completes the proof.
\end{proof}

\begin{Lemma}\label{ch-poly}
	The characteristic polynomial of $B_n$, the $\beta$-matrix of order $n$, is 
	\begin{equation*}
		p_n(t)=\sum_{j=0}^{n}t^j-\sum_{i=1}^{n}\sum_{j=0}^{n-i}t^{i+j-1}\beta^{-i}.
	\end{equation*}
\end{Lemma}
\begin{proof}
	In order to compute $p_n(t)$, we use basic properties of the determinant and then the Laplace expansion together with the auxiliary result obtained in Lemma \ref{aux-determinant}. $p_n(t)$ is the determinant of the matrix $tI_n-B_n$, that is
	\begin{equation*}
	\begin{pmatrix}
	t-\beta^{-1}+1&-\beta^{-1}+1&-\beta^{-1}+1&-\beta^{-1}+1&\cdots&-\beta^{-1}+1\\
	-\beta^{-2}-1&t-\beta^{-2}&-\beta^{-2}&-\beta^{-2}&\cdots&-\beta^{-2}\\
	-\beta^{-3}&-\beta^{-3}-1&t-\beta^{-3}&-\beta^{-3}&\cdots&-\beta^{-3}\\
	-\beta^{-4}&-\beta^{-4}&-\beta^{-4}-1&t-\beta^{-4}&\ddots&-\beta^{-4}\\
	\vdots&\ddots&\ddots&\ddots&\ddots&\vdots\\
	-\beta^{-n}&\cdots&\cdots&-\beta^{-n}&-\beta^{-n}-1&t-\beta^{-n}
	\end{pmatrix}.
	\end{equation*}
	First, we subtract the last column from all the previous ones, which does not change the determinant of the matrix. As a consequence we obtain 
		\begin{equation*}
	A_n=\begin{pmatrix}
	t&0&\cdots&\cdots&\cdots&0&-\beta^{-1}+1\\
	-1&t&0&\cdots&\cdots&0&-\beta^{-2}\\
	0&-1&t&0&\cdots&0&-\beta^{-3}\\
	\vdots&\ddots&\ddots&\ddots&\ddots&\vdots&\vdots\\
	0&\cdots&0&-1&t&0&-\beta^{-n+2}\\
	0&\cdots&0&0&-1&t&-\beta^{-n+1}\\
	-t&\cdots&\cdots&\cdots&-t&-1-t&t-\beta^{-n}
	\end{pmatrix}.
	\end{equation*}
	In order to expand using the Laplace rule along the last column, we have to study the cofactors of type $(i,n)$ of $A_n$. Actually, if we remove the last column and the i-th row from the matrix $A_n$, we then obtain a lower triangular block matrix, with two blocks on the diagonal. The first one is a matrix of order $i-1$, it is itself lower triangular with $t$ on the diagonal entries, while the second one is $M_{n-i}$, i.e. the matrix of Lemma \ref{aux-determinant} of order $n-i$. For instance, the minor $(3,6)$ of $A_6$ is given by the matrix below
	\begin{equation*}
		\begin{pmatrix}
		t&0&0&0&0\\
		-1&t&0&0&0\\
		0&0&-1&t&0\\
		0&0&0&-1&t\\
		-t&-t&-t&-t&-1-t\\
		\end{pmatrix}.	
	\end{equation*}
	Therefore, the determinant of the minor $(i,n)$ of $A_n$ is the product of the determinants of the two diagonal blocks, that is 
	\begin{equation*}
		{\rm det}((A_n)_{i,n})=t^{i-1}{\rm det}(M_{n-i})=(-1)^{n-i}t^{i-1}\sum_{j=0}^{n-i}t^j=(-1)^{n-i}\sum_{j=0}^{n-i}t^{i+j-1}.
	\end{equation*}
	Finally, we can then compute $p_n(t)$ as
	\begin{eqnarray*}
		{\rm det}(A_n)=\sum_{i=1}^{n}(-1)^{n+i}a_{i,n}{\rm det}((A_n)_{i,n})=(-1)^{n+1}(1-\beta^{-1})(-1)^{n-1}\sum_{j=0}^{n-1}t^j+\\
		+\sum_{i=2}^{n-1}(-1)^{n+i}(-\beta^{-i})(-1)^{n-i}\sum_{j=0}^{n-i}t^{i+j-1}+(-1)^{2n}(t-\beta^{-n})t^{n-1}=\\
		=\sum_{j=0}^{n-1}t^j+t^n+\sum_{i=1}^{n}(-\beta^{-i})\sum_{j=0}^{n-i}t^{i+j-1}=\sum_{j=0}^{n}t^j-\sum_{i=1}^{n}\sum_{j=0}^{n-i}t^{i+j-1}\beta^{-i},
	\end{eqnarray*}
	where we write the final expression of $p_n(t)$, by separating the terms with powers of $\beta$ from the others.
\end{proof}

\begin{Observation}\label{Obs-invertible}
	The constant term of $p_n(t)$ is $1-\beta^{-1}$, the difference between the constant terms of the two sums, and since $\beta>1$ we conclude that, for every n, $B_n$ does not have $0$ as an eigenvalue and, as a consequence, is invertible.
\end{Observation}
	
	The sequence of characteristic polynomials $\{p_n(t)\}_{n\in\N}$ has a very particular structure and we want to exploit it in order to derive information on the distribution of the roots.
	For every $n$, we define the following polynomials
	\begin{eqnarray*}
		q_n(t)=\sum_{j=0}^{n}t^j,\\
		r_n(t)=\sum_{i=1}^{n}\sum_{j=0}^{n-i}t^{i+j-1}\beta^{-i}.
	\end{eqnarray*}
	It immediately follows that
	\begin{equation*}
	p_n(t)=q_n(t)-r_n(t).
	\end{equation*}

\begin{Lemma}\label{inside-circle}
	
	For every $t\in\mathbb{C}$ of absolute value less than $1$, we have
	\begin{eqnarray*}
		\lim_{n \to +\infty} q_n(t)=\frac{1}{1-t},\\
		\lim_{n \to +\infty} r_n(t)=\frac{1}{1-t}\frac{1}{\beta-t},
	\end{eqnarray*}
	and the convergence is uniform in every compact set $K$ contained in the open ball $B(0,1)\subset\mathbb{C}$.
\end{Lemma}
\begin{proof}
	The result of $q_n(t)$ follows trivially from the basic properties of power series and by the fact that $q_n(t)$ is the partial sum of order $n$ of the geometric series, which has a radius of convergence equal to $1$. Similarly, we can also note that $r_n(t)$ is related to the partial sum of the Cauchy product of two geometric series, of ratio $t$ and $t/\beta$. By rearranging, we can write
	\begin{equation*}
		r_n(t)=\sum_{i=1}^{n}\sum_{j=0}^{n-i}t^{i+j-1}\beta^{-i}=\frac{1}{\beta}\sum_{i=0}^{n-1}\sum_{j=0}^{n-i-1}t^{i+j}\beta^{-i}=
		\frac{1}{\beta}\sum_{i=0}^{n-1}\sum_{j=0}^{n-i-1}t^{j}\bigg(\frac{t}{\beta}\bigg)^i,
	\end{equation*}
	where the last double sum is the sum of all the possible products of powers of $t$ and $t/\beta$ such that the exponents add up to at most $n-1$. For this reason, apart from the factor $1/\beta$, $r_n(t)$ is exactly the partial sum of order $n-1$ of the Cauchy product of the two geometric series, as required. As far as uniform convergence is concerned, the radius of convergence of the geometric series of ratio $t$ is $1$, whilst for the one of ratio $t/\beta$ it is $\beta$. By Mertens' theorem on the convergence of Cauchy products, the radius of convergence of the product series is at least the minimum between $1$ and $\beta$, which is $1$. Therefore, for every complex $t$ of modulus less than $1$, we find
	\begin{equation*}
		\lim_{n \to +\infty} r_n(t)=\frac{1}{\beta}\frac{1}{1-t}\frac{1}{1-t/\beta}=\frac{1}{1-t}\frac{1}{\beta-t}
	\end{equation*}
	and the convergence is uniform in every compact set $K$ contained in the open disk $B(0,1)\subset\mathbb{C}$. This completes the proof.
\end{proof}
In order to deduce the complete result on the clustering of the eigenvalues, we would also like to study the exterior of the unitary circle, which is not possible by a direct limit process. Nevertheless, given a polynomial $p(t)$ of degree $n$ without $0$ as a root, the function $t^np(1/t)$ is again a polynomial of degree exactly $n$, whose roots are the reciprocals of the roots of $p(t)$. This transformation swaps the inside and the outside of the unit circle. In this way, we may study what happens to the roots outside the unit circle, for $n$ large enough, by studying what happens in the limit. 
Thus we define
\begin{eqnarray*}
	\tilde{p}_n(t)=t^np_n(1/t),\\
	\tilde{q}_n(t)=t^nq_n(1/t)=q_n(t),\\
	\tilde{r}_n(t)=t^nr_n(1/t).\\
\end{eqnarray*}
As in the previous case, we are interested in computing the limit.
\begin{Lemma}\label{outside-circle}
	For every $t\in\mathbb{C}$ of absolute value less than $1$, we find
	\begin{equation*}
		\lim_{n \to +\infty} \tilde{r}_n(t)=\frac{t}{1-t}\frac{1}{\beta-1}
	\end{equation*}
	and the convergence is uniform in every compact set $K$ contained in the open disk $B(0,1)\subset\mathbb{C}$.
\end{Lemma}
\begin{proof}
	As in the previous case, the result is trivial from the same properties of power series and Cauchy product series, once we have rearranged $\tilde{r}_n(t)$. Therefore, following similar steps as in Lemma \ref{inside-circle}, we observe
	\begin{equation*}
		\tilde{r}_n(t)=t^nr_n(1/t)=\sum_{i=1}^{n}\sum_{j=0}^{n-i}t^{n+1-i-j}\beta^{-i}=\frac{1}{\beta}
		\sum_{i=0}^{n-1}\sum_{j=0}^{n-i-1}t^{n-i-j}\beta^{-i}.
	\end{equation*}
	Then we switch the index $j$ in the internal sum with the index $k=n-i-j-1$ which allows us to go through that sum in the inverse order
	\begin{equation*}
	\frac{1}{\beta}\sum_{i=0}^{n-1}\sum_{j=0}^{n-i-1}t^{n-i-j}\beta^{-i}=\frac{1}{\beta}\sum_{i=0}^{n-1}\sum_{k=0}^{n-i-1}t^{k+1}\beta^{-i}=\frac{t}{\beta}\sum_{i=0}^{n-1}\sum_{k=0}^{n-i-1}t^{k}\beta^{-i}.
	\end{equation*}
	As with the previous case, the last sum is the sum of all the possible products of powers of $t$ and $1/\beta$ such that the exponents add up to at most $n-1$. For this reason, apart from the factor $t/\beta$, the polynomial $\tilde{r}_n(t)$ is the partial sum of the Cauchy product series of two geometric series of ratio $t$ and $1/\beta$. By Mertens' theorem, for every $t\in\mathbb{C}$ of absolute value less than $1$, we deduce
	\begin{equation*}
	\lim_{n \to +\infty} \tilde{r}_n(t)=\frac{t}{\beta}\frac{1}{1-t}\frac{1}{1-1/\beta}=\frac{t}{1-t}\frac{1}{\beta-1}
	\end{equation*}
		and the convergence is uniform in every compact set $K$ contained in the open disk $B(0,1)\subset\mathbb{C}$. 
With the latter the proof is concluded.
\end{proof}

The last element for the proof of the main result is Hurwitz' Theorem, which establishes relations between the zeros of a sequence of holomorphic functions and the zeros of the uniform limit.
\begin{Theorem}\label{Hurwitz}
	Let $\{f_n\}_{n\in\N}$ be a sequence of holomorphic functions on a connected open set $A$ that converge uniformly on compact subsets of $A$ to a holomorphic function $f$, which is not constantly zero on $A$. If $f$ has a zero of order $m$ at $z_0$ then for every small enough $\epsilon > 0$ and for sufficiently large $n\in\N$, $f_n$ has precisely $m$ zeros in the disk defined by $\abs{z-z_0}<\epsilon$, counted with their multiplicity. Furthermore, these $m$ zeros converge to $z_0$ as $n\to\infty$.
\end{Theorem}
With all the instruments developed so far, we are ready to prove our main result, that is Theorem \ref{thm-strong-cluster}.
\begin{proof}
	Given the $\beta$ matrix-sequence $\{B_n\}_{n\in\N}$, we computed the sequence of characteristic polynomials $p_n(t)$. By Lemma \ref{inside-circle}, for every $t\in\mathbb{C}$ of modulus less than $1$, we have that
	\begin{equation*}
		\lim_{n \to +\infty} p_n =\lim_{n \to +\infty} q_n(t) - r_n(t)=\frac{1}{1-t}
		 -\frac{1}{1-t}\frac{1}{\beta-t}=\frac{1}{1-t}\frac{\beta-1-t}{\beta-t}
	\end{equation*}
	and the convergence is uniform in every compact $K$ subset of the open ball $B(0,1)$. Similarly, by Lemma \ref{outside-circle}, for every $t\in\mathbb{C}$ of modulus less than $1$ we obtain
	 \begin{equation*}
	 \lim_{n \to +\infty} \tilde{p}_n =\lim_{n \to +\infty} \tilde{q}_n(t) - \tilde{r}_n(t)=\frac{1}{1-t}
	 -\frac{t}{1-t}\frac{1}{\beta-1}=\frac{1}{1-t}\frac{\beta-1-t}{\beta-1}
	 \end{equation*}
	 and the convergence is once again uniform, when we restrict to a compact subset in the open ball $B(0,1)$.
	 We call $p(t)$ and $\tilde{p}(t)$ the two limit functions, namely
	 \begin{eqnarray*}
	 	p(t)=\frac{1}{1-t}\frac{\beta-1-t}{\beta-t},\\
	 	\tilde{p}(t)=\frac{1}{1-t}\frac{\beta-1-t}{\beta-1}.
	 \end{eqnarray*}
 We observe that	$p$ and $\tilde{p}$ are meromorphic with a single zero equal to $\beta-1$. Therefore, if we take $\beta\in\left[2,+\infty\right)$, and $\epsilon>0$ the sequences $p_n(t)$ and $\tilde{p}_n(t)$ converges uniformly to $p(t)$ and $\tilde{p}(t)$ in the compact closed disk centered at $0$ of radius $1-\epsilon$ and the two limit functions never take the value $0$ in the ball. Thus, there exists $N_{\epsilon}\in\N$ such that, for every $n\geq N_{\epsilon}$, neither $p_n(t)$ nor $\tilde{p}_n(t)$ has roots in the compact set. This completes the proof for the case $\beta\in\left[2,+\infty\right)$ since for those values of $\beta$ all the eigenvalues of $\{B_n\}_{n\in\N}$ and their reciprocals are strongly clustered at $\{s\in\mathbb{C}:|s|\geq1\}$, with no outliers, for every $\epsilon >0$.\\
 	Now, let $\beta\in\left(1,2\right)$. We cannot use the same argument itself, but we need a certain refinement of it, since the two limit functions $p(t)$ and $\tilde{p}(t)$ have a unique zero in the unit open disc, namely $\tilde{z} = \beta-1$. By derivation of the two limit functions we obtain
 	\begin{eqnarray*}
 		p'(t)=\frac{1}{(1-t)^2}\frac{(\beta-t)^2-(\beta-t)-(1-t)}{(\beta-t)^2},\\
 		\tilde{p}'(t)=\frac{1}{(1-t)^2}\frac{(\beta-1)-t-(1-t)}{(\beta-1)}.
 	\end{eqnarray*}
 	Thus
 	\begin{eqnarray*}
 		p'(\tilde{z})=\frac{\beta-2}{\beta^2}\neq0,\\
 		\tilde{p}'(\tilde{z})=\frac{\beta-2}{\beta^2(\beta-1)}\neq0.
 	\end{eqnarray*}
 Consequently, $p(t)$ and $\tilde{p}(t)$ have a unique single zero inside the open unit disk. By Theorem \ref{Hurwitz}, for sufficiently large $n\in\N$, $p_n(t)$ and $\tilde{p}_n(t)$ have exactly one single root in the open disk centered at $\tilde{z}$ of radius $\epsilon$ and those converge to $\tilde{z}$, when $n\to\infty$. In addition, for $n$ large enough, the functions $p_n(t)$ and $\tilde{p}_n(t)$ cannot have roots in any compact subset of the open unit disk not containing $\tilde{z}$ due to uniform convergence. Hence, we conclude that all but two eigenvalues must lie, at some point, in an $\epsilon$-neighborhood of the unit circle.
 	In order to complete the proof, we show that the polynomials, for $n$ large enough, have a real positive root not approaching $S^1$. By the previous argument, such a root must approach the critical value $\tilde{z}$.
 	First, by evaluation, we have
 	\begin{eqnarray*}
 		p_n(0)=\frac{\beta-1}{\beta}>0,\\
 		\tilde{p}_n(0)=1>0.
 	\end{eqnarray*} 
 	Secondly, working on the limit functions as real functions, we obtain
 	\begin{eqnarray*}
 		\lim_{t\to1^-}p(t)=\lim_{t\to1^-}\frac{1}{1-t}\frac{\beta-1-t}{\beta-t}=-\infty,\\
 		\lim_{t\to1^-}	\tilde{p}(t)=\lim_{t\to1^-}\frac{1}{1-t}\frac{\beta-1-t}{\beta-1}=-\infty.
 	\end{eqnarray*}
 	If we consider the polynomials as real functions, we have uniform convergence to the limit functions on $\left[0,\delta\right]$ for every $0<\delta<1$. Therefore, taking $\delta$ such that $\tilde{z}<\delta$, and noting that the polynomials have a negative sign for $n$ large enough when evaluated at $\delta$, once again by uniform convergence. For such $n$, we can use Bolzano's theorem on the zeros of continuous functions to ensure that both $p_n(t)$ and $\tilde{p}_n(t)$ have a zero in $\left[0,\delta\right]$. This completes the proof since $\delta$ is arbitrary but fixed and so these zeros cannot approach $S^1$ for $n\to\infty$. As a conclusion, by all previous arguments, they approach $\tilde{z} = \beta-1$. This completes the case of $\beta\in\left(1,2\right)$ and the proof of the theorem is concluded.
\end{proof}

\section{Numerics, Concluding Remarks, Open Problems}\label{end}

In this final section, we will present some numerical evidence in support of Theorem \ref{thm-strong-cluster}. Then we give additional related results and draw conclusions, as well as stressing a few open questions.

\begin{figure}[ht]
    \centering
    \begin{minipage}{0.49\textwidth}
        \centering
        \includegraphics[width=0.99\textwidth]{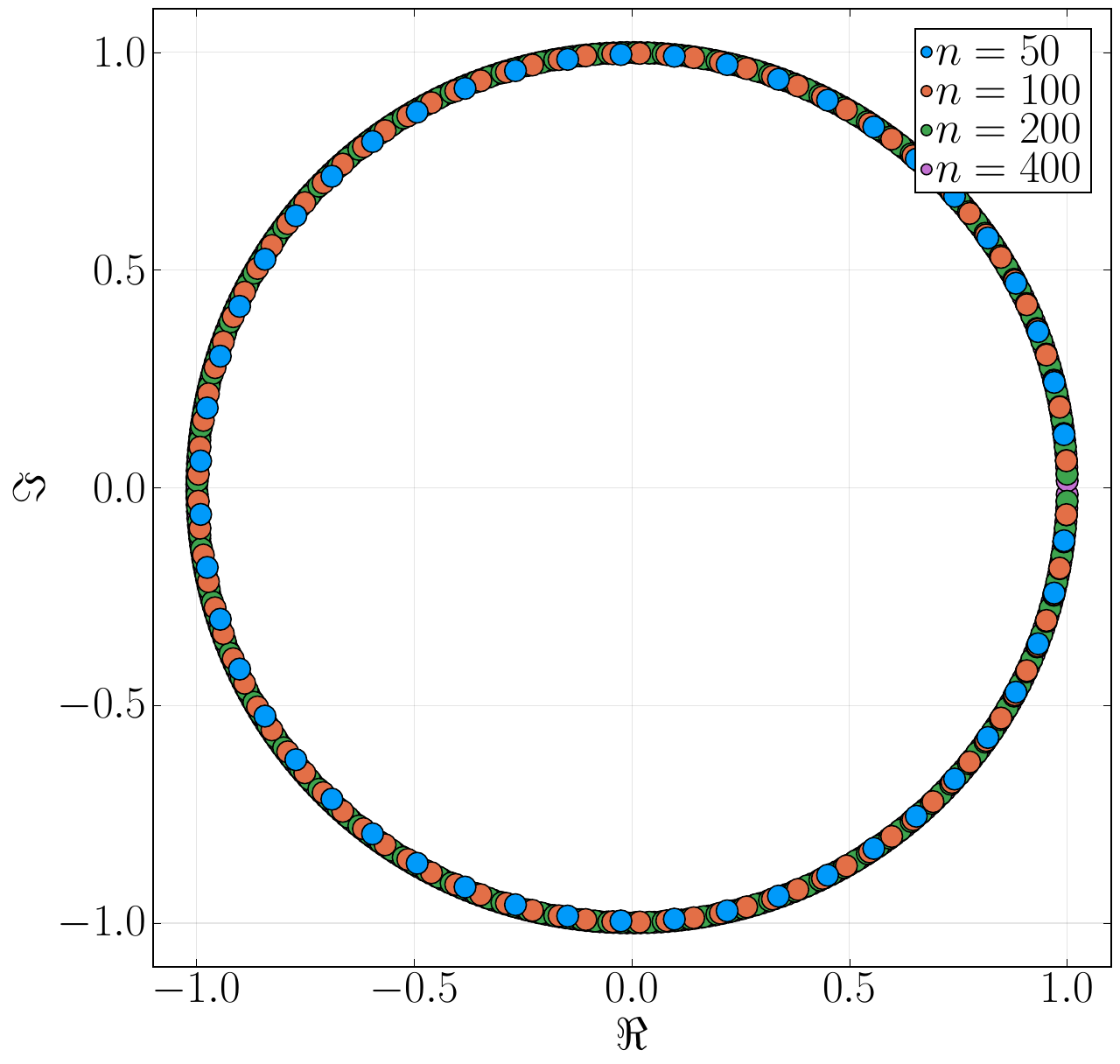} 
        \caption{$\beta=5$}\label{fig1-david}
    \end{minipage}\hfill
    \begin{minipage}{0.49\textwidth}
        \centering
        \includegraphics[width=0.99\textwidth]{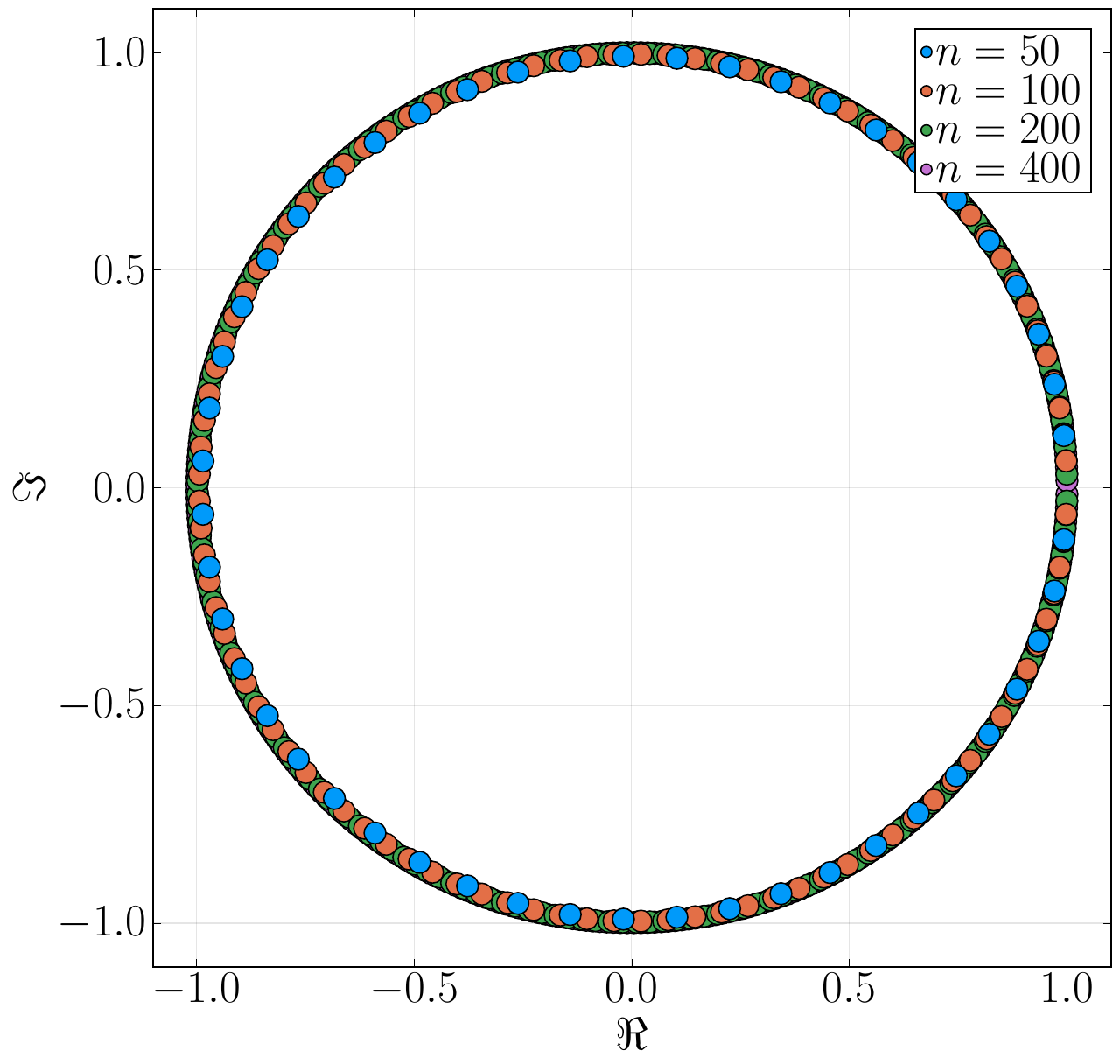} 
        \caption{$\beta=3$}\label{fig2-david}
    \end{minipage}
\end{figure}

\begin{figure}[ht]
	\centering
	\includegraphics[width=\textwidth]{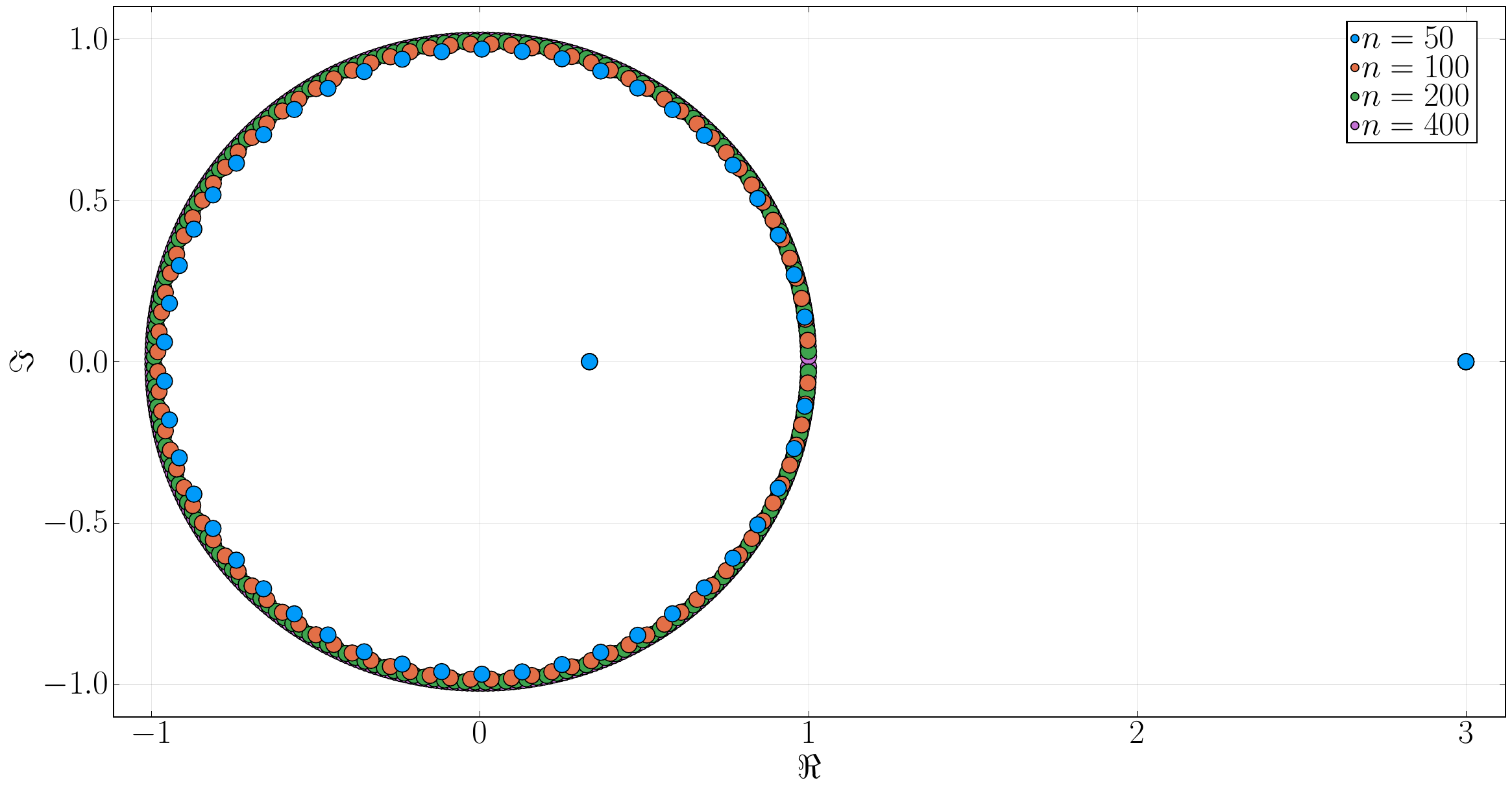}
	\caption{$\beta =\frac{4}{3}$}
	\label{fig3-david}
\end{figure}

Firstly, we plot the spectra $B_n$ for different values of $\beta$ and for different sizes $n$. The computations are performed with high precision due to the non-normal character of the matrices.  The following cases are considered:
\begin{enumerate}
\item plot of all the eigenvalues for $\beta=5$, $n=50, 100, 200, 400$;
\item plot of all the eigenvalues for $\beta=3$, $n=50, 100, 200, 400$;
\item plot of all the eigenvalues for $\beta=4/3$, $n=50, 100, 200, 400$.
\end{enumerate}
\ 
\noindent 
The results are visualized in Figure \ref{fig1-david}, Figure \ref{fig2-david}, and Figure \ref{fig3-david}.

In the first two figures with $\beta=5\ge 2$ and $\beta=3\ge 2$, according to Theorem \ref{thm-strong-cluster}, no ouliers are expected and the eigenvalues strongly cluster at the complex unit circle. As can be seen, the numerical experiments strongly agree with the previous theoretical results even for small values of $n$. In the third figure we have taken $\beta=4/3\in (1,2)$ so that exactly two outliers are expected, converging to $\beta-1=1/3$ and $(\beta-1)^{-1}=3$, respectively. Also in this case, the numerics confirm the theoretical findings and we observe the limit behavior already for moderate sizes. In fact, for the third experiment and for the largest outlier $\lambda_{M}$ we have the following impressive results:
\, 
\\ \,
\begin{itemize}
\item $n=50$, $\lambda_{M}=2.99999796124162120902813536126303334491749260835507
$;
\item $n=100$, $\lambda_{M}=
2.9999999999988454072132625253185082984139093876636
$;
\item $n=200$, $\lambda_{M}=
2.9999999999999999999999996296987491150278175529157
$;
\item $n=400$, $\lambda_{M}=
2.99999999999999999999999999999999999999999999999996
$.
\end{itemize}
\, \\
\\ \,
Finally, one can immediately see that the proof of Theorem \ref{thm-strong-cluster} does not apply to the limit case where $\beta=1$. However, it suggests that one outlier tends to $0$  and the other to $+\infty$, when $n$ tends to infinity. For this we go back to the matrix $B_n$, taking $\beta=1$, which takes the form

\begin{equation*}
	 	B_n=\begin{pmatrix}
	 	0 & 0 & 0& 0 &\cdots& 0\\
	 	2&1 & 1& 1 &\cdots& 1\\
	 	1 & 2& 1 & 1&\cdots& 1\\
	 	1 & 1& 2& 1&\ddots& 1\\
	 	\vdots&\ddots&\ddots&\ddots&\ddots&\vdots\\
	 	1 &\cdots&\cdots& 1 & 2 & 1
	 	\end{pmatrix}.
	 	\end{equation*}
The matrix is in lower block triangular form with zeros in the first row. Thus the minimal outlier is clearly equal to $\beta-1=0$ with multiplicity $1$ since the other block diagonal submatrix is invertible and takes the form $X_{n-1}=T_{n-1}(e^{\mathbf{i}\theta})+ e e^{T}$, with $e=e^{(n-1)}$ being the vector of size $n-1$ with all ones. The related eigenvectors $w$ are such that $B_n w=0$ and can be expressed as
\[
w=(w)_1 
\begin{pmatrix}
1 \\
s
\end{pmatrix}, \ \ \ \ s=-X_{n-1}^{-1}r, \ \ \ r=e+e_1,
\]
with $(w)_1\neq 0$ the first component of $w$.

Given the lower block triangular structure of $B_n$, the other eigenvalues, including the largest outlier, must be found in the matrix $X_{n-1}$.
We first observe that $X_{n-1}$ is a matrix with all positive entries. Therefore, by Perron-Fronenius theory (see \cite{P-F} and references therein), its spectral radius coincides with the largest eigenvalue, i.e. the largest outlier, which is unique. Since the spectral radius is reached by a single eigenvalue with algebraic multiplicity equal to $1$, the power method is convergent (see \cite{eig-comput}). We then use this fact for estimating the largest eigenvalue which takes the expression
\begin{equation}\label{asympt:exp}
\lambda_{M}=n-1/n+c_2/n^2 +c_3/n^3+ \cdots
\end{equation}
where $c_0=0$, $c_1=-1$, and the other $c_j$, $j\in \N$, are constants independent of $n$. In fact, we refer to Table \ref{table-david}, where we report the first $5$ iterates of the power method applied to $X_{n-1}$ starting with initial vector $v_0=e$, $v_{k+1}=X_{n-1}v_k$, $k\ge 0$, and where $r_k = \frac{(v_{k+1})_1}{(v_k)_1}$, $k\ge 0$, is the scalar quantity converging to the maximal eigenvalue. We notice that $r_k$ is always well defined since all the involved vectors $v_k$, $k\ge 0$, have strictly positive entries. In conclusion, by observing that $c_1=-1$ is negative, we see $\lambda_M<n$. The fact that $\lambda_M<n$ is in theoretical accordance with the third Gerschgorin Theorem (see \cite{Va04}), taking into account that $X_{n-1}$ is irreducible and its first disk is contained in the open part of the others. Furthermore, $\lambda_M$ tends to $+\infty$ as $n$ tends to $\infty$, in agreement with the qualitative suggestion given by Theorem \ref{thm-strong-cluster}. However, we have even more information, since $\lambda_M=\lambda_M(n)$ has an asymptotic expansion which can be employed for using the linear in time matrix-less methods introduced by Ekstr\"om and coauthors (see \cite{SE2,SE1} for the seminal papers). In this direction see Table \ref{table-davidBIS} which is a numerical counterpart of the symbolic computations in Table \ref{table-david}.

\begin{table}[ht]\centering
\begin{tabular}{|l|l|l|}
\hline
$k$ & $(v_k)_1$   & $r_k = \frac{(v_{k+1})_1}{(v_k)_1}$ \\ \hline
1 & $n-1$   &  $n - \frac{1}{n-1}$    \\ \hline
2 & $n^2 -n-1$ &  $n - \frac{n}{n^2 -n-1}$  \\ \hline
3 & $n^3 -n^2 -2n$ &  $n - \frac{n-1}{n(n-2)}$  \\ \hline
4 & $n^4 -n^3 -3n^2 +1$ & $n-\frac{n^3 + 2n + 1}{n^4 -n^3 -3n^2 +1}$ \\ \hline
5 & $n^5 -n^4 - 4n^3 +3n+1$   &  $n -\frac{n^3+n^2+2n}{n^4 -4n^2 - 4n -1}$  \\ \hline
\end{tabular}\caption{Sequence showing the approximate largest eigenvalue}\label{table-david}
\end{table}

\begin{table}[ht]\centering
\begin{tabular}{|l|l|l|}
\hline
$n$ & $0=c_0\approx \lambda_M(n)-n$   & $-1=c_1\approx  n(\lambda_M(n)-n)$ \\ \hline
50  & $-0.0204166702$   &  $-1.0208335106$    \\ \hline
100 & $-0.0101020409$   &  $-1.0102040921$  \\ \hline
200 & $-0.0050252525$   &  $-1.0050505056$  \\ \hline
400 & $-0.0025062814$   &  $-1.0025125628$ \\ \hline
\end{tabular}\caption{Sequence showing numerically the asymptotic expansion in (\ref{asympt:exp})}\label{table-davidBIS}
\end{table}

To conclude, in the current paper, we have studied the eigenvalue behaviour for a class of non-normal matrix-sequences as defined in \ref{Beta-mat-seq}, which represent a one-rank correction depending on a parameter $\beta>1$ of the basic Toeplitz sequence generated by $e^{\mathbf{i}\theta}$. We have proven a strong clustering at the range of the generating function $e^{\mathbf{i}\theta}$ i.e. at the complex unit circle for every fixed $\beta>1$.  
For $\beta\ge 2$ there are no outliers, while for $\beta \in (1,2)$ only two outliers are present, which are both real, positive and have a finite limit equal to $\beta-1$ and $(\beta-1)^{-1}$, respectively. The conditioning measured in spectral norm grows at least as 
\[
\left[\max\{\beta-1,(\beta-1)^{-1}\}\right]^2
\]
for any choice of $\beta>1$.
Numerical experiments have confirmed the theoretical findings with high convergence speed, i.e. already for quite moderate matrix orders, and also in the limit case of $\beta=1$.

As already observed in the Introduction, the problem looks mathematically innocent, but indeed is quite challenging since all the sophisticated asymptotic linear algebra machinery for deducing the weak clustering is not easy to apply in the current setting and at most we may hope for weak clustering results. 
In our derivations, we have determined a representation of the characteristic polynomials, which have been worked out with the help of several tools from complex analysis. 

Further distribution results in the Weyl sense can be obtained by joining old findings in the theory of distribution of the zeros of polynomial sequences and potential theory \cite{saff,saff-book}, as well as newer research in asymptotic linear algebra (see \cite{book GLT I} and references therein). We resume then in Theorem \ref{Weyl-distributions} and in the subsequent two remarks.

\begin{Theorem}\label{Weyl-distributions}
Let $\{B_n\}_{n\in\N}$ be the $\beta$ matrix-sequence given in Definition \ref{Beta-mat-seq}.
\begin{itemize}
\item Assume that $|\beta|\ge 1$. Then $\{B_n\}_{n\in\N}$ is distributed in the eigenvalue sense as $e^{\mathbf{i}\theta}$ on $[-\pi,\pi]$ i.e.
for any continuous function $F$ defined on the complex field and with bounded support we have
\[
    \lim_{n\to\infty}\frac1{n}\sum_{i=1}^{n}F(\lambda_i(B_n))=\frac1{2\pi}\int_{-\pi}^\pi F(e^{\mathbf{i}\theta}){\rm d}\theta. 
\]
\item Assume that $\beta$ is any nonzero complex number. Then $\{B_n\}_{n\in\N}$ is distributed in the singular value sense as $1$ on $[-\pi,\pi]$
i.e. for any continuous function $F$ defined on the non-negative real numbers and with bounded support we have
\[
    \lim_{n\to\infty}\frac1{n}\sum_{i=1}^{n}F(\sigma_i(B_n))=\frac1{2\pi}\int_{-\pi}^\pi F(1){\rm d}\theta =F(1). 
\]
\end{itemize}
\end{Theorem}

\begin{proof}
By Lemma \ref{ch-poly}, the characteristic polynomial of $B_n$ is 
\[		p_n(t)=\sum_{j=0}^{n}t^j-\sum_{i=1}^{n}\sum_{j=0}^{n-i}t^{i+j-1}\beta^{-i}
\]
so that, as the leading coefficient is $1$ and for any complex $\beta$ with $|\beta|\ge 1$, we have
\[
1 \le \|p_n\|_{L^\infty(S^1)} \le n+n^2,
\] 
$S^1$ being the unit circle in $\mathbb{C}$. Hence 
\[
\lim_{n\to\infty}  \|p_n\|_{L^\infty(S^1)}^{1/n}= 1
\]
and the first part of the theorem is then proved by using results from classical potential theory (see \cite{saff-book} and references therein).

The second part has an elementary proof due to the interlacing theorems of the singular values.  In fact, looking at Definition \ref{Beta-mat-seq}, the singular values of $T_n(e^{\mathbf{i}\theta})$ are given by $1$ repeated $n-1$ times and by $0$. Additionally, $B_n=T_n(e^{\mathbf{i}\theta})+(v-e_1)e^{T}$ is a rank one correction of $T_n(e^{\mathbf{i}\theta})$ for any $\beta\neq 0$, so that all the singular values are equal to $1$ except for two possible outliers.
A direct computation allows us to conclude that 
\[
\lim_{n\to\infty}\frac1{n}\sum_{i=1}^{n}F(\sigma_i(B_n))=F(1).
\]

\end{proof}

\begin{Observation}\label{Obs-Weyl1}
	When a matrix-sequence is distributed as a function in the eigenvalue sense, there is a corresponding weak eigenvalue clustering at the essential range of the function.   When a matrix-sequence is distributed as a function in the singular value sense, there a corresponding weak singular value clustering at the essential range of the modulus of the function. In this sense Theorem \ref{Weyl-distributions} implies a weak clustering at $1$ of the singular values of $\{B_n\}_{n\in\N}$ and a weak clustering on the complex unit circle of the eigenvalues, since the complex unit circle is the range of the function $e^{\mathbf{i}\theta}$. However, via distribution results it is difficult to obtain the strong clustering results obtained in the present work. Conversely, for the singular values the strong clustering is for free for any $\beta\neq 0$ thanks to the interlacing theorems for the singular values as already indicated in the proof of the second part of Theorem \ref{Weyl-distributions}.
\end{Observation}

\begin{Observation}\label{Obs-Weyl2}
	No matter how the parameter $\beta$ is chosen the matrix $B_n$ is always non-normal. By the Schur normal form, a characterization of normality is that the moduli of the eigenvalues are equal to the singular values. Hence, by Theorem \ref{Weyl-distributions}, when $|\beta|\ge 1$, what is not true for the matrix $B_n$ is true in an asymptotic sense for the matrix-sequence $\{B_n\}_{n\in\N}$: in fact for $|\beta|\ge 1$ the sequences of sets $\{\{|\lambda_i(B_n)|\}_{i=1}^2\}_{n\in\N}$ and $\{\{\sigma_i(B_n)\}_{i=1}^2\}_{n\in\N}$ are both distributed as the constant function $1$ and the difference sequence $\{\{\sigma_i(B_n)-|\lambda_i(B_n)|\}_{i=1}^2\}_{n\in\N}$ is distributed as the constant function $0$. 
Hence we can state that $\{B_n\}_{n\in\N}$ is a quasi-normal matrix-sequence.
\end{Observation} 

Other generalizations can be considered. If instead of $B_n=T_n(e^{\mathbf{i}\theta})+(v-e_1)e^{T}$ in (\ref{Bn compact}) we consider $\hat B_n=T_n(e^{\mathbf{i}\theta})-e_1e^{T} +v z^T$, where $z$ in any vector with zeros and ones different from $e$, then the correction is of rank two and the complexity grows. Some cases may be recovered using the same analysis proposed in the present work, while the general setting is still an open question.
A further direction to investigate is that suggested by relation (\ref{asympt:exp}) and Table \ref{table-davidBIS}, regarding almost closed forms for the eigenvalues of $B_n$ with respect to $\beta$ and $n$: this would open the door to linear in time fast computations of the entire spectrum of $B_n$ for very large $n$.

\section*{Acknowledgments}
The authors would like to thank Giovanni Canestrari for introducing us to the literature on $\beta$ maps, and Giovanni Canestrari, Sven-Erik Ekstr\"om, and Carlo Garoni for fruitful discussions on various aspects of the related matrix-theoretic problem. 
The Referees helped us in improving the quality of the presentation and stimulated us in making the content stronger (Theorem \ref{Weyl-distributions} and the related discussion in Observation \ref{Obs-Weyl1} and Observation \ref{Obs-Weyl2}): hence we warmly thank them  for their careful reading and helpful suggestions.\\
The work of Stefano Serra-Capizzano was supported in part by INDAM-GNCS and was funded from the European High-Performance Computing Joint Undertaking  (JU) under grant agreement No 955701. The JU receives support from the European Union’s Horizon 2020 research and innovation programme and Belgium, France, Germany, Switzerland. Furthermore Stefano Serra-Capizzano is grateful for the support of the Laboratory of Theory, Economics and Systems – Department of Computer Science at Athens University of Economics and Business.

\end{document}